\documentclass[aos,MSNbibl,seceqn,dvips]{arximspdf}
\usepackage{graphicx}

\doi{10.1214/11-AOS927}
\volume{39}
\issue{5}
\pubyear{2011}
\firstpage{2244}
\lastpage{2265}



\begin{document}
\begin{frontmatter}

\title{Erich Leo Lehmann---A glimpse into his life and~work}
\runtitle{Erich Leo Lehmann---A glimpse into his life and work}

\begin{aug}
\author[A]{\fnms{Javier} \snm{Rojo}\corref{}\thanksref{t1}\ead[label=e1]{jrojo@rice.edu}}
\runauthor{J. Rojo}
\affiliation{Rice University}
\address[A]{Department of Statistics\\
Rice University\\
Houston, Texas 77005\\
USA\\
\printead{e1}} 
\end{aug}

\thankstext{t1}{Supported in part by the
NSF Grants DMS-10-19634, DMS-08-51368 and by the National Security
Agency through Grant H9823000601-0099.}

\received{\smonth{9} \syear{2011}}

%
\begin{abstract}
Through the use of a system-building approach, an approach that
includes finding common ground for the various philosophical paradigms
within statistics, Erich L. Lehmann is responsible for much of the
synthesis of classical statistical knowledge that developed from the
Neyman--Pearson--Wald school. A biographical sketch and a brief summary
of some of his many contributions are presented here. His complete
bibliography is also included and the references present many other
sources of information on his life and his work.
\end{abstract}

%

\end{frontmatter}

\section{Biographical sketch}
Erich L. Lehmann was born in Strasbourg on November 20th, 1917. He
passed away in Berkeley, California on the morning of September
12th, 2009. His family left Germany in 1933, as the Nazis came to
power, to settle in Switzerland. He spent five years in Z\"{u}rich and
two years in Trinity College in Cambridge studying mathematics. Under
the United States French immigration quota---Strasbourg was, by then,
part of France as a consequence of the Versailles treaty---he arrived
in New York at the end of 1940. Edmund Landau, the famous number
theorist, was an acquaintance of the Lehmann family and had suggested
Trinity College as the place Erich should go to study mathematics.
Landau died in 1938 from a heart attack, but his wife wrote a letter of
introduction for Erich to take to Landau's G\"{o}ttingen colleague
Richard Courant who was now in New York developing what became the
Courant Institute. Courant, having offered the option to ``live in New
York or in the United States,'' and Erich having opted for the latter,
recommended the University of California as an up-and-coming good
place. Erich arrived in Berkeley, California in January 1, 1941.

Erich's first order of business was to speak with Griffith C. Evans,
chair of the mathematics department, who immediately accepted him as a
probationary graduate student. The probationary status resulted from
Erich not having a degree. Evans, who had been recruited from the
mathematics department at Rice Institute---now Rice University---had a
broad vision for mathematics and had the intention of hiring Ronald A.
Fisher, whom he knew. However, a visit by Fisher to Berkeley did not go
well. The news of Jerzy Neyman's successful visit to the United States,
culminating with a set of lectures at the U.S. Department of
Agriculture, reached Evans who in 1937 offered Neyman a job in the
mathematics department at the University of California without having
met him. With the advent of the second world war, Evans advised Erich
that it might be a good idea to move from mathematics to some other
area---perhaps physics or statistics---that could be more useful to
the war efforts. Erich, not being fond of physics, opted for
statistics. His initial experiences, however, led him to second-guess
his decision. In Lehmann (2008B), Erich writes that ``statistics did
not possess the beauty that I had found in the integers and later in
other parts of mathematics. Instead, ad hoc methods were used to solve
problems that were messy and that were based on questionable
assumptions that seemed quite arbitrary.'' (Hereafter, a reference
followed by ``B'' indicates book reference in Section 9; a bracketed
reference [x] refers to that numbered reference in Section 8; other
references appear at the end of this work.) After some soul-searching,
he decided to go back to mathematics and approached the great logician
Alfred Tarski. Tarski accepted him as a student, but before Erich had
an opportunity to let Evans and Neyman know about his decision, Neyman
offered him a job as a lecturer with some implicit potential for the
position to become permanent. Feeling that this represented a great
opportunity to become part of a community, something that Erich very
much desired at that point in time, he decided to take the offer and
abandoned his plans for returning to mathematics. In 1942 Erich
received an M.A. degree in mathematics, and was a teaching assistant in
the Statistical Laboratory from 1942 to 1944 and from 1945 to 1946.

These early years as a graduate student and a teaching member of the
depart\-ment---while sharing office space with Charles Stein, Joseph
Hodges and Evelyn Fix---helped to forge lifetime friendships and
productive collaborations. After he spent the year from August 1944 to
August 1945 stationed in Guam as an operations analyst in the United
States Air Force, Erich returned to Berkeley and started working on a
thesis problem proposed by Pao-Lu Hsu in consultation with Neyman. The
problem was in probability theory---some aspect of the moment
problem---and after obtaining some results and getting ready to write
them up,
Erich discovered that his results were already in Markov's work. The
situation became complicated as Neyman was invited to supervise the
Greek elections. Before leaving, Neyman asked Hsu if he could provide
another thesis topic for Erich. Hsu obliged but was not able to
supervise Erich's thesis, as he followed Hotelling from Columbia to
North Carolina and then decided to go back to China. Neyman turned to
George P\`{o}lya at Stanford for help. Weekly meetings with P\`{o}lya,
commuting between Berkeley and Stanford, finally yielded a thesis.
Meanwhile, Neyman was back from Greece after being relieved of his
duties for insubordination. Neyman had felt that the elections were
rigged and decided to check by himself. When asked to stop, he refused.
This turn of events allowed Neyman to be back in Berkeley for Erich's
examination. Thus, in June of 1946, Erich obtained his Ph.D. degree
with a thesis titled ``Optimum tests of a certain class of hypotheses
specifying the value of a correlation coefficient.''

Erich was not the first of Neyman's Berkeley Ph.D. students, but he was
the first one to be hired by the mathematics department. He held the
title of assistant professor of mathematics from 1947 to 1950, and
spent the first half of 1950--1951 as a visiting associate professor at
Columbia, and as a lecturer at Princeton during the second half of that
year. Partly to allow more time for the tumultuous situation created in
Berkeley by the anti-Communist loyalty oath to settle down, and partly
to make a decision on an offer from Stanford, Erich spent the year of
1951--1952 as a visiting associate professor at Stanford. Erich decided
to go back to Berkeley, but not before he was able to persuade Neyman
not to require him to do consulting work for the statistical
laboratory. (Stanford's offer explicitly mentioned that Erich was not
expected to do any applied work.) On his return to Berkeley in 1952,
Erich was promoted to associate professor of mathematics, and then in
1954 was promoted to professor of mathematics. In 1955, after Evans
stepped down as chair of mathematics, thus providing Neyman with his
opportunity for a new department of statistics, Erich's title changed
to professor of statistics. In 1988, Erich became professor emeritus
and then from 1995 to 1997 he was distinguished research scientist at
the Educational Testing Service (ETS). In spite of his retirement in
1988, Erich continued to be professionally active and a regular
participant in the social life of the department. Despite offers from
Stanford in 1951 and from the Eidgen\"ossische Technische Hochschule
(ETH) in 1959, and except for short stints at Columbia, Princeton,
Stanford and ETS, Erich lived in Berkeley
from his arrival on January 1st, 1941 until his death on September
12th, 2009.

\section{Honors, awards, service to the profession}
Erich Lehmann's towering contributions to statistics have received many
well-deserved accolades. Erich was an elected fellow of the Institute
of Mathematical Statistics (IMS) and of the American Statistical
Association (ASA), and he was an elected member of the International
Statistical Institute. Remarkably, he was the recipient of three
Guggenheim Fellowships (1955, 1966 and 1980) and two Miller Institute
for Basic Research Professorships (1962 and 1972). The IMS honored him
as the Wald lecturer in 1964---the title of his lectures being ``Topics
in Nonparametric Statistics.'' This was followed in 1988 by the
Committee of Presidents of Statistical Societies (COPSS) R. A. Fisher
Memorial Lecture entitled ``Model Specification: Fisher's views and
some later strategies.'' In 1975 Erich was elected fellow of the
American Academy of Arts and Sciences and in 1978 he was elected member
to the National Academy of Sciences. Election as an Honorary Fellow of
the Royal Statistical Society followed in 1986 and the ASA recognized
him with the Wilks Memorial Award in 1996. His life-long work was
recognized with two Doctorates \textit{honoris causa}, the first from
the University of Leiden in 1985, and the second from the University of
Chicago in 1991. The honor from Leiden carries with it the distinction
of being the first Dr. \textit{h. c.} granted by the University of
Leiden to a mathematician in a century, the previous one having been
awarded to Stieltjes in 1884. In 1997, to celebrate Erich's 80th
birthday, the Berkeley statistics department instituted the Lehmann
fund to provide support for students. In 2000 Erich became the first
Goffried Noether Award recipient and lecturer for his influential work
in Nonparametrics. His Noether lecture, entitled ``Parametrics Versus
Nonparametrics: Two Alternative Methodologies,'' formed the basis for
an invited paper with discussion in the \textit{Journal of
Nonparametrics}
(JNPS) in 2009 [121]. Posthumously, Erich received the best JNPS paper
award for 2009. His students and colleagues honored him with a set of
reminiscences in 1972 (J. Rojo, ed.), a \textit{Festschrift for Erich
L. Lehmann} organized by Bickel, Doksum and Hodges in 1982 [see also
Bickel, Doksum and Hodges (1983)], and a series of \textit{Lehmann
Symposia}, organized by Rojo and Perez-Abreu in 2002, and Rojo in 2004,
2007 and 2011. Perhaps surprisingly, although he was honored with the
Fisher lecture, he never received the honor of being the Neyman
lecturer. It may be surmised that Erich's lack of affinity for applied
work impeded his being so honored.

Erich served the profession well. Although initially reluctant to serve
as chair of the statistics department at Berkeley, he did so from
1973--1976. And he did it very well. Brillinger (2010) writes:

\begin{quote}
He had always refused previously for a variety of reasons. He did
it so well that I sometimes thought that he must have thought through
how a Chair should behave and put his conclusions into practice. For
example, to the delight of visitors and others he was in the coffee
room each day at 10 a.m. He focused on the whole department---staff,
students, colleagues and visitors.
\end{quote}

During 1960--1961, Erich was IMS President and was a leader in the
internationalization of the IMS [see, e.g., Lehmann (2008B) and van
Zwet (2011)]. He was a member of the Executive committee of the Miller
Institute (1966--1970), and a member of the committee of visitors to the
Harvard Department of Statistics (1974--1980) and Princeton (1975--1980).
He served as Editor of the \textit{Annals of Mathematical Statistics} from
1953--1955 and as Associate Editor from 1955--1968. He was invited to
stay on for a second term as Editor but, after accepting, had to
decline. For details see Lehmann (2008B) and van Zwet (2011).


\section{Books and their translations}

In his youth, Erich Leo Lehmann had a desire to become a writer. In
Lehmann (2008B), he wrote, ``\textit{My passion was German literature, my
dream to become a writer, perhaps another Thomas Mann or Gottfried Keller}.''
Surely it was this passion that drove Erich to write his successful and
influential books. The list includes:

\begin{enumerate}
\item\textbf{Testing Statistical Hypotheses}. Three editions (1959, 1986,
2005). The 2005 edition is joint with Joseph P. Romano. The 1959
edition was translated into Russian (1964), Polish (1968) and Japanese.
\item\textbf{Basic Concepts of Probability and Statistics}, with Joseph
L. Hodges. Two editions (1964, 1970). Reprinted in 2005 as part of the
SIAM series Classics in Applied Mathematics. The book was translated
into Hebrew (1972), Farsi (1994), Italian (1971) and Danish (1969).

\item\textbf{Elements of Finite Probability}, with Joseph L. Hodges. Two
editions (1965, 1970).

\item\textbf{Nonparametrics: Statistical Methods Based on Ranks}, with
the assistance of H. J. M. D'Abrera. Hardcover edition (1975) by
Holden-Day. Paperback edition (1998) by Prentice-Hall, Inc., Simon \&
Schuster, and then by Springer Science in 2006. The book was translated
into Japanese (1998).
\item\textbf{Theory of Point Estimation}. Two editions (1983,
1998---with George Casella). The 1983 edition was translated into Russian
(1991), and the 1998 edition into Chinese (2004).
\item\textbf{Elements of Large-Sample Theory}, 1999.
\item\textbf{Reminiscences of a Statistician: The Company I Kept}, 2008.
\end{enumerate}

Additionally, Erich collaborated with Judith M. Tanur on the book
\textit{Statistics}: \textit{A~Guide to the Unknown}. This book went
through several editions and translations [Chinese (1980) and Spanish
(1992)]. Spin-offs from this book were two other books with similar
titles: \textit{Statistics}: \textit{A Guide to the Study of the
Biological and Health Sciences} and \textit{Statistics}: \textit{A
Guide to Political and Social Issues}, both published in 1977, and on
which Erich collaborated. Erich served as co-editor or special editor.
The complete list of books and their translations is given in
Section~9.

\begin{figure}
\begin{tabular}{@{}cc@{}}
\multicolumn{2}{@{}c@{}}{
\includegraphics{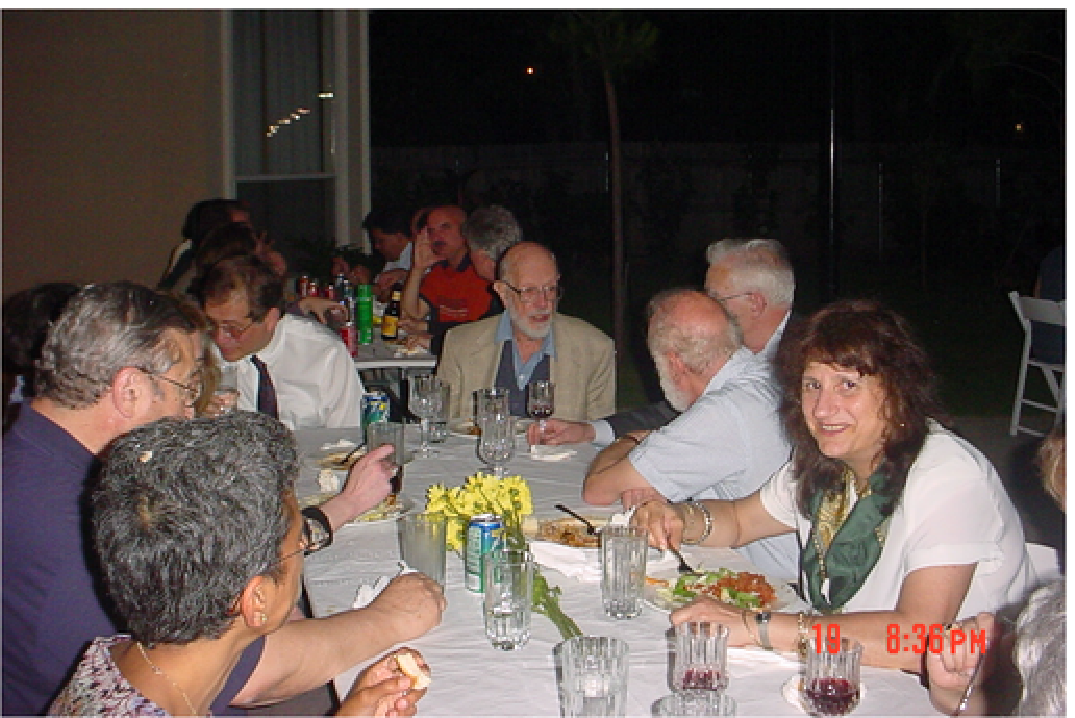}
}\\
\multicolumn{2}{@{}p{\columnwidth}@{}}{\fontsize{9pt}{11pt}\selectfont{\textit{Erich L. Lehmann
at dinner during the 2nd Lehmann
Symposium---May 2004. Shown also in the picture}, \textit{David Cox}, \textit{Ingram
Olkin}, \textit{Peter Bickel}, \textit{Shulamith Gross}, \textit{Emmanuel Parzen and Loki Natarajan.
Kjell Doksum}, \textit{Joseph Romano and Gabriel Huerta are seen in the
background.}}}
\end{tabular}
\\[6pt]
\begin{tabular}{@{}c@{\qquad}c@{}}
\multicolumn{1}{@{}c}{
\includegraphics{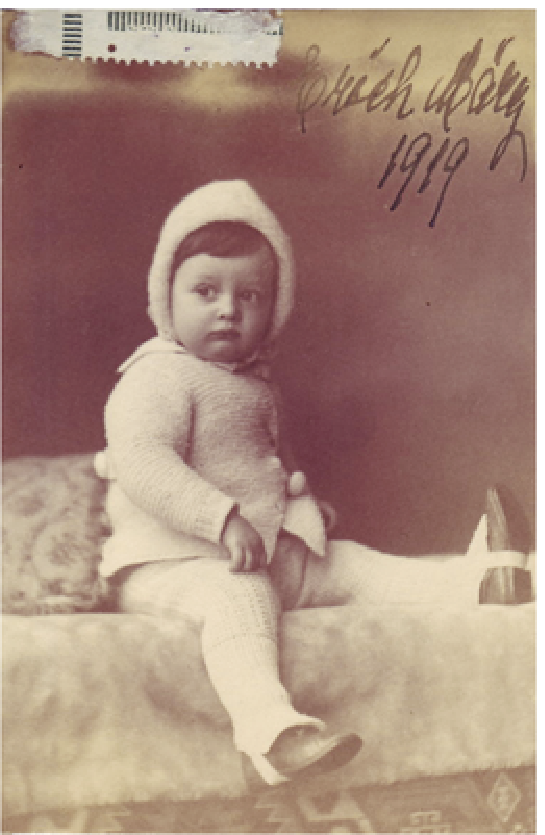}
} &
\multicolumn{1}{c@{}}{
\includegraphics{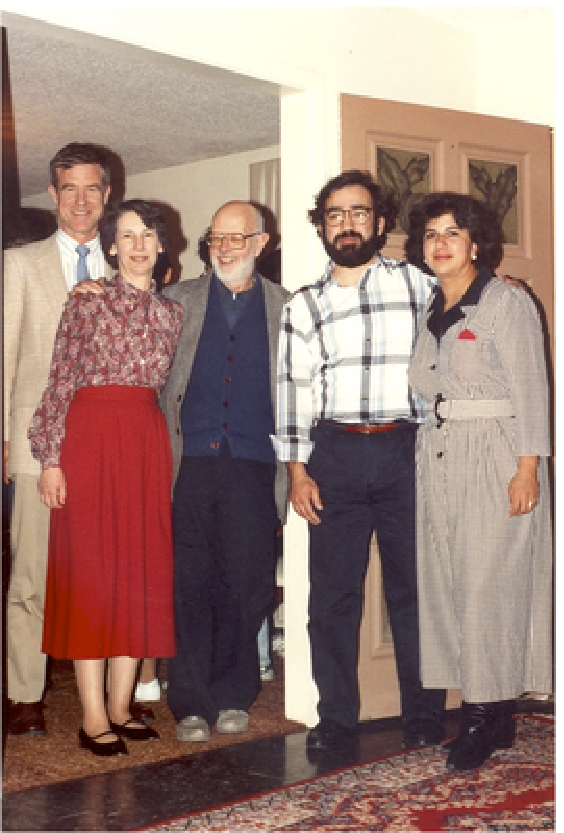}
}\\
\multicolumn{1}{c}{\fontsize{9pt}{11pt}\selectfont{\textit{Erich L. Lehmann in 1919.}}} &
\multicolumn{1}{p{160pt}@{}}{\fontsize{9pt}{11pt}\selectfont{\textit{Erich L. Lehmann in 1992 in El Paso},
\textit{TX.
Also in the picture Juliet Shaffer}, \textit{Javier and Ma
Luisa Rojo}, \textit{and Simon Bernau.}}}
\end{tabular}
\\[10pt]
\textbf{Erich L. Lehmann in 1919, 1992 and 2004.}
\end{figure}

The book \textit{Fisher}, \textit{Neyman}, \textit{and the Creation of
Classical Statistics} has now been published posthumously by Springer,
Lehmann (2011B). Erich was finishing the manuscript at the time of his
death. Juliet Shaffer worked diligently after Erich's passing to bring
the book to publication form. Fritz Scholz continues work on a
revision, started before Erich's death, of the \textit{Nonparametrics}:
\textit{Statistical Methods Based on Ranks} book. The revision
incorporates the use of R and the book is expected to be completed in
two years.

\section{Technical work}
Erich's contributions are multifaceted and too many to do justice to in
the allotted space. A more extensive
and careful assessment of his work is provided in Rojo (2011). Here,
only a small part of his work will be briefly
reviewed. Some of his ground-breaking work in nonparametric statistics
is discussed in this issue by van Zwet (2011).

\subsection{Early work}
While still a graduate student at Berkeley, Erich submitted a paper
that was published in 1947 [2], in which the issue of what to do
when a uniformly most powerful (UMP) test does not exist is discussed.
Erich proposed that, due to the many tests available to choose from,
one must reduce attention
to a class of tests $\mathcal{F}$ with the property that for any test
$\phi$ not in $\mathcal{F}$, there is a test $\phi^*$ in $\mathcal{F}$
with a power function at least as good as that of $\phi$. And if $\phi
_1$ and $\phi_2$ are two tests in $\mathcal{F}$, then neither one
dominates the other. In addition, the paper characterizes the class
$\mathcal{F}$ for a special case. Erich recognized that the class
$\mathcal{F}$ may still be too large to offer much relief in finding a
good solution and, therefore, other information or principles may be
needed to further narrow down the class $\mathcal{F}$. Thus, the
concept of minimal complete classes, that plays a fundamental role in
the theory of statistical decisions of Wald (1950), was born in this paper.

In his book \textit{Statistical Decision Functions} (1950), Wald credits Lehmann:

\begin{quote}
The concept of complete class of decision functions was
introduced by Lehmann, and the first result regarding such classes is
due to him [30]$\ldots.$
\end{quote}

Interestingly, Neyman was not impressed by this work. In DeGroot (1986)
Erich states:

\begin{quote}
I wrote it
up---it was just a few pages---and said to Neyman that
I would like to publish it. He essentially said, ``It's
junk. Do not bother.'' But I sent it in to Wilks anyway.
\end{quote}

Some of Erich's early work was motivated by the work of Hsu (1941)
that dealt with optimal properties of the likelihood ratio test in the
context of analysis of variance. In Lehmann (1959) [34], Erich shows
that these optimal properties are consequences of the fact that the
test is uniformly most powerful invariant. In addition, the paper
unified optimality results of Kiefer (1958) for symmetrical
nonrandomized designs, and optimality results
of Wald (1942) for the analysis of variance test for the general
univariate linear hypothesis.
Hsu (1941) also proposed a method for finding all similar tests.
Lehmann (1947) [3] extended Hsu's results to the composite null
hypothesis problem, and ideas in Hsu (1941) motivated the concept of
completeness in Lehmann and Scheff{\'e} (1950) [12]. Lehmann and
Scheff{\'e} (1950) [12] and Lehmann and Scheff{\'e} (1955) [26]
provided a comprehensive study of the concepts of similar regions and
sufficient statistics. Together with Lehmann and Stein (1950) [11],
where uniformly minimum variance unbiased estimators are discussed in
the sequential sampling context, these papers provide the final word on
certain problems in hypotheses testing and estimation.

\subsection{Minimaxity and admissibility}

Hodges and Lehmann (1950, 1951, 1952) [10, 13, 16], provided minimax
estimators for several examples and the admissibility of minimax
estimators and connections with Bayes estimators were discussed. In
Hodges and Lehmann (1950) [10], a minimax estimator for the probability
of success $p$ in a binomial experiment is obtained by considering the
Bayes estimator with respect to a beta conjugate prior that yields a
Bayes estimator with constant risk. The minimax estimator thus found is
admissible due to the uniqueness of the Bayes estimator. The results
are extended to the case of two independent binomial distributions and
a minimax estimator is obtained for the difference of the probability
of successes when the sample sizes are equal. The question of whether a
minimax estimator exists for the difference of the success
probabilities for unequal sample sizes remains an open problem. The
papers also consider the nonparametric case, and methods for deriving
nonparametric minimax estimators are provided under certain conditions.
The concept of complete classes having been formalized by Wald (1950),
the paper also shows that, for convex loss functions, the class of
nonrandomized estimators is essentially complete.
Hodges and Lehmann (1951) [13] used a different approach to obtain
minimax and admissible estimators when the loss function is a weighted
squared error loss. The method requires the solution of a differential
inequality involving the lower bound for the Mean Squared Error.
Various sequential problems were discussed and minimax estimators were derived.
Hodges and Lehmann (1952) [16] proposed finding estimators whose
maximum risk does not exceed the minimax risk by more than a given
amount $r$. Under this restriction it was proposed to find the \textit
{restricted Bayes solution} with respect to some prior distribution
$\lambda$. That is, find $\delta_0$ that minimizes $\int R(\delta,\theta
)\,d\lambda(\theta)$ subject to $\sup_{\theta}R(\delta,\theta)\le r$.
Conditions were discussed for the existence of restricted Bayes
estimators and several examples were provided that illustrate the
method. It was argued that Wald's theory can be extended to obtain
results for these restricted Bayes procedures.

Wald (1950) obtained the existence of least favorable distributions
under the assumption of a compact parameter space. Lehmann (1952) [17]
addressed this issue and, in the case of hypothesis testing and, more
generally, in the case where only a finite number of decisions are
available, Lehmann weakened the conditions for the existence of least
favorable distributions.
Lehmann and Stein (1953) [20] proved the admissibility of the most
powerful invariant test when testing certain hypotheses in the location
parameter family context.

\subsection{Hypothesis testing} Erich's work on hypothesis testing is
well known. Here some aspects of that work are briefly reviewed.

\subsubsection{Composite null hypotheses} Lehmann (1947) [3] and
Lehmann and Stein (1948) [5] studied the problem of testing a composite
(null) hypothesis. The 1947 paper extends the work of Scheff{\'e}
(1942). Suppose that $\Theta$ is a $k$-dimensional parameter space. Let
$\Theta_0$ be the subset of $\Theta$ given by $\{\bar{\theta}\in\Theta
\dvtx\theta_i=\theta_i^0\}$, for one $i=1,\ldots,k$. Then the null
hypothesis $H_0\dvtx\bar{\theta}\in\Theta_0$ is an example of a composite
(null) hypothesis with one constraint, and the parameters $\theta_j, j
\ne i$, are nuisance parameters. Neyman (1935) provided Type B regions
for the case of a single nuisance parameter. These results were
extended by Scheff{\'e} to the case of several nuisance parameters
(under $H_0$), and Scheff{\'e} provided sufficient conditions for these
Type B regions to also be Type B$_1$ (uniformly most powerful unbiased)
regions. Lehmann (1947) [3] utilized Neyman and Pearson's (1933) and
Hsu's (1945) methods to determine the totality of similar regions and
extended Scheff{\'e}'s results to obtain uniformly most powerful tests
against one-sided alternatives. Hsu's method was also employed to
obtain UMP regions in cases, for example, location and scale
exponential and uniform distributions, where Neyman and Pearson's
method does not apply. The above approach is not as fruitful in the
case of more than one constraint, but results of Hsu (1945) are useful
in this regard.

In Lehmann and Stein (1948) [5] the problem of testing a composite
hypothesis against a single alternative is addressed by relaxing the
condition of similarity to one requiring only that $\int_{\Omega
^*}f(x)\,d\le\alpha$ or all $f\in\mathcal{F}$, where $\Omega^*$ denotes
the critical region of the test. Adapting the Neyman--Pearson lemma to
hold in this case, sufficient conditions for the existence of most
powerful tests were derived. The results for Student's problem, with
composite null hypothesis given by the normal family with mean $0$ and
unknown variance, and the simple alternative hypothesis given by the
normal distribution with known parameters were somewhat surprising; see
Lehmann (2008B), page 48.

\subsubsection{Likelihood ratio tests}

Lehmann (1950, 1959, 2006) [9, 34, 118] deal with
the likelihood ratio principle for testing. Although
this principle is ``intuitive'' and provides ``reasonable'' tests,
it is well known that it may fail. The papers
examine different aspects of the problem focusing on the optimality of
the likelihood ratio test in some cases, and in its total failure in
other cases.

Lehmann (1959) [34]
considered a class of invariant tests endowed with
an order that satisfies certain properties. It was then shown that, in this
case, the likelihood ratio test's optimality properties follow
directly from the fact
that the test is uniformly most powerful invariant. See also Section 4.1.

In Lehmann (2006) [118] and Lehmann (1950) [9], properties of tests
produced by other approaches
are examined and compared to the likelihood ratio tests. For example,
when the testing
problem remains invariant with respect to a transitive
group of transformations,
the \textit{likelihood averaged or integrated with respect to an
invariant measure approach} in
Lehmann (2006) [118] produces tests that
turn out to be uniformly at least as powerful as
the corresponding likelihood ratio test, with the former being
strictly better except when
the two coincide; and in the absence of invariance, the
proposed approach continues to improve on the likelihood ratio test
for many
cases. Lehmann (1950) [9] was discussed in Section 4.1.

\section{Orderings of probability distributions}
Lehmann's work on orderings of probability distributions was motivated
in part from the need to study properties of power functions. Thus,
Lehmann (1955) [27] discussed the stochastic and monotone likelihood
ratio orderings. The latter plays a fundamental role in the theory of
uniformly most powerful tests and both can be characterized in terms of
the function $K(u)=GF^{-1}(u)$; see, for example, Lehmann and Rojo
(1992) [98]. It is this function $K$ that also plays a fundamental role
in the Lehmann Alternatives and, hence, is also connected with the Cox
proportional hazards model and has now spilled over to the literature
on Receiving Operating Characteristic (ROC) curves.
A different collection of partial orderings between distributions $F$
and $G$ can be defined in terms of the function $K^*=F^{-1}G(x)$.
Bickel and Lehmann (1979) [64] considered the dispersive ordering
defined by requiring that $K^*(y)-K^*(x)\le y-x$ for all $y>x$, and
considered several of its characterizations. This concept is
equivalent, under some conditions, to a tail-ordering introduced by
Doksum (1969). This function, $K^*$, is also useful in comparing
location experiments (Lehmann (1988) [85]).

Lehmann (1966) [47] introduced concepts of dependence for random
variables $(X,Y)$. This work has attracted a lot of attention in the
literature from applied probabilists and statisticians alike.

\section{Philosophical work}

Erich believed in the frequentist interpretation of probability and in
the Neyman--Pearson--Wald school of optimality, but recognized that both
perspectives have their limitations. See, for example, page 188 of
Lehmann (2008B). Bickel and Lehmann (2001) [110, 111] discussed some of
the philosophical shortcomings of a frequentist interpretation of
probability. Erich felt that optimality considerations achieve
solutions that may lack robustness and other desirable properties.
His work on foundational issues focused on the following: (i)~model
selection; (ii) frequentist statistical inference; (iii) Bayesian
statistical inference; and (iv) exploratory data analysis.

Restricting attention to (ii), (iii) and (iv), Erich viewed the
trichotomy as being ordered by the level of model assumptions made.
Thus, (iv) is free of any model assumptions and allows the data to
speak for itself, while the frequentist approach relies on a
probability model to evaluate the procedures under consideration. The
Bayesian approach, in addition, brings in the prior distribution. Erich
felt that none of these approaches is perfect. Motivated by this state
of affairs, Lehmann (1985, 1995) [82, 104] developed ideas that bridge
the divide created by the heated philosophical debates. Lehmann (1985)
[82] discussed how the Neyman--Pearson--Wald approach contributes to the
exploration of underlying data structure and its relation with Bayesian
inference. Lehmann (1995) [104] continued with this line of thought:

\begin{quote}
In practice, the three approaches can often fruitfully interact, with
each benefiting from considerations of the other points of view. It
seems clear that model-free data analysis, frequentist and Bayesian
model-based inference and decision making each has its place. The
question appears not to be---as it is often phrased---which is the
correct approach but in what circumstances each is most appropriate.
\end{quote}

Erich's balanced view of foundational issues is appealing. His work
reflects the belief that no single paradigm is totally satisfactory.
Rather than exacerbating their differences through heated debates, he
proposed that a fruitful approach is possible by consolidating the good
ideas from (ii), (iii) and (iv)---with (iii) serving as a bridge that
connects all three.
Although his original position was solidly in the frequentist camp, he
shifted, somewhat influenced by classical Bayesian ideas. However, he
felt that a connection with the radical Bayesian position was more
challenging. He states in Lehmann (1995) [104] that ``\textit{bridge
building to the \textup{``}radical\textup{'' [}Bayesian\textup{]}
position is more difficult}.'' A
definition of the radical Bayesian position is not provided, but it can
be surmised that this refers to a paradigm that insists on the
elicitation of a prior distribution at all costs. In Lehmann (2008B),
he writes:

\begin{quote}
However, it seems to me that the strength of these beliefs tends
to be rather fuzzy, and not sufficiently well defined and stable to
assign a definite numerical value to it. If, with considerable effort,
such a value is elicited, it is about as trustworthy as a confession
extracted through torture.
\end{quote}

\section{Ph.D. students} I first attended U.C. Berkeley during the
Fall of 1978. My first course was statistics 210 A---the first quarter
of theoretical statistics. The recollections of my days as a student
during that first quarter, followed by two more quarters of theoretical
statistics---statistics 210 B and C---all taught by Erich, are very
vivid. During
that first academic year, I was very impressed with Erich's lecturing
style. He would present the material without unnecessarily dwelling too
long on technical details, and in such a way that connections with
previous material seemed virtually seamless. It was quite enjoyable to
follow ``the story'' behind the theory. His lectures were so perfectly
organized even when only using a few notes on his characteristic
folded-in-the-coat's-pocket-yellow sheets! Regarding teaching, Erich wrote
in Lehmann (2008B):

\begin{quote}
While I eschewed very large courses, I loved the teaching that occurred
at the other end of the spectrum. Working on a one-on-one basis with
Ph.D. students was, for me, the most enjoyable and rewarding aspect of
teaching. At the same time, it was an extension of my research, since
these students would help me explore areas in which I was working at
the time.
\end{quote}

This love for one-on-one teaching produced a total of 43 Ph.D.
students. Curiously, two
of Erich's Ph.D. students obtained their degrees from Columbia rather
than from Berkeley. That these students graduated from Columbia, rather
than from Berkeley, resulted from a confluence of circumstances.
Although Erich had received an invitation from Wald to visit Columbia
during the 1949--1950 academic year, Erich had to postpone his visit to
Columbia for the following year since Neyman took a sabbatical during
the 1949--1950 academic year. After Wald's tragic and untimely death,
two of Wald's students approached Erich with a request to become his
students. These students are marked with an asterisk in the following
table that presents the names and dissertation titles, by year of
degree, for all 43 of Erich's Ph.D. students.

\begin{enumerate}[1950]
\item[1950]\textbf{Colin Ross Blyth}\\
\textit{I. Contribution to the Statistical Theory of the Geiger--Muller
Counter}; \\
\textit{II. On Minimax Statistical
Decision Procedures and Their Admissibility.}

\item[1953]\textbf{Fred Charles Andrews}\\
\textit{Asymptotic Behavior of Some Rank Tests for Analysis of Variance.}
\\
\textbf{Jack Laderman*}\\
\textit{On Statistical Decision Functions for Selecting One of k Populations.}

\item[1954]\textbf{Hendrik Salomom Konijn}\\
\textit{On the Power of Some Tests for Independence.}

\item[1955]\textbf{Allan Birnbaum*}\\
\textit{Characterizations of Complete Classes of Tests of Some
Multiparametric Hypotheses, with
Applications to Likelihood Ratio Tests.}\\
\textbf{Balkrishna V. Sukhatme}\\
\textit{Testing the Hypothesis that Two Populations Differ Only in Location.}

\item[1959]\textbf{V. J. Chacko }\\
\textit{Testing Homogeneity Against Ordered Alternatives.}

\item[1961]\textbf{Piotr Witold Mikulski}\\
\textit{Some Problems in the Asymptotic Theory of Testing Statistical
Hypotheses.}

\item[1962]\textbf{Madan Lal Puri}\\
\textit{Asymptotic Efficiency of a Class of C-Sample Tests.}
\\
\textbf{Krishen Lal Mehra}\\
\textit{Rank Tests for Incomplete Block Designs. Paired-Comparison Case.}
\\
\textbf{Subha Bhuchongkul Sutchritpongsa }\\
\textit{Class of Non-Parametric Tests for Independence in Bivariate
Populations.}
\\
\textbf{Shishirkumar Shreedhar Jogdeo}\\
\textit{Nonparametric Tests for Regression Models.}

\item[1963]\textbf{Peter J. Bickel}\\
\textit{Asymptotically Nonparametric Statistical Inference in the
Multivariate Case.}
\\
\textbf{Arnljot H\o yland }\\
\textit{Some Problems in Robust Point Estimation.}

\item[1964]\textbf{Milan Kumar Gupta}\\
\textit{An Asymptotically Nonparametric Test of Symmetry.}
\\
\textbf{Madabhushi Raghavachari}\\
\textit{The Two-Sample Scale Problem When Locations are Unknown.}
\\
\textbf{Ponnapalli Venkata Ramachandramurty}\\
\textit{On Some Nonparametric Estimates and Tests in the
Behrens--Fisher Situation.}
\\
\textbf{Vida Greenberg}\\
\textit{Robust Inference in Some Experimental Designs.}

\item[1965]\textbf{Kjell Andreas Doksum}\\
\textit{Asymptotically Minimax Distribution-Free Procedures.}
\\
\textbf{William Harvey Lawton}\\
\textit{Concentration of Random Quotients.}

\item[1966]\textbf{Shulamith Gross}\\
\textit{Nonparametric Tests When Nuisance Parameters Are Present.}
\\
\textbf{Bruce Hoadley}\\
\textit{The Theory of Large Deviations with Statistical Applications.}
\\
\textbf{Gouri Kanta Bhattacharyya}\\
\textit{Multivariate Two-Sample Normal Scores Test for Shift.}
\\
\textbf{James Nwoye Adichie}\\
\textit{Nonparametric Inference in Linear Regression.}
\\
\textbf{Dattaprabhakar V. Gokhale}\\
\textit{Some Problems in Independence and Dependence.}

\item[1968]\textbf{Frank Rudolph Hampel}\\
\textit{Contributions to the Theory of Robust Estimation.}

\item[1969]\textbf{Wilhelmine von Turk Stefansky}\\
\textit{On the Rejection of Outliers by Maximum Normed Residual.}
\\
\textbf{Neil H. Timm}---Co-advisors Erich Leo Lehman and Leonard
Marascuilo\\
\textit{Estimating Variance--Covariance and Correlation Matrices from
Incomplete Data.}
\\
\textbf{Louis Jaeckel}\\
\textit{Robust Estimates of Location.}

\item[1971]\textbf{Friedrich Wilhelm Scholz}\\
\textit{Comparison of Optimal Location Estimators.}
\\
\textbf{Dan Anbar}\\
\textit{On Optimal Estimation Methods Using Stochastic Approximation
Procedures.}

\item[1972]\textbf{Michael Denis Stuart}\\
\textit{Components of 2 for Testing Normality Against Certain Restricted
Alternatives.}
\\
\textbf{Claude L. Guillier}\\
\textit{Asymptotic Relative Efficiencies of Rank Tests for Trend Alternatives.}
\\
\textbf{Sherali Mavjibhai Makani}\\
\textit{Admissibility of Linear Functions for Estimating Sums and
Differences of Exponential Parameters.}

\item[1973]\textbf{Howard J. M. D'Abrera}\\
\textit{Rank Tests for Ordered Alternatives.}

\item[1974]\textbf{Hyun-Ju Yoo Jin}\\
\textit{Robust Measures of Shift.}

\item[1977]\textbf{Amy Poon Davis}\\
\textit{Robust Measures of Association.}

\item[1978]\textbf{Jan F. Bjornstad}\\
\textit{On Optimal Subset Selection Procedures.}

\item[1981]\textbf{William Paul Carmichael}\\
\textit{The Rate of Weak Convergence of a Vector of U-Statistics Generated
by a Single Sample.}
\\
\textbf{David Draper}\\
\textit{Rank-Based Robust Analysis of Linear Models.}

\item[1982]\textbf{Wei-Yin Loh}\\
\textit{Tail-Orderings on Symmetric Distributions with Statistical
Applications.}

\item[1983]\textbf{Marc J. Sobel}\\
\textit{Admissibility in Exponential Families.}

\item[1984]\textbf{Javier Rojo}\\
\textit{On Lehmann's General Concept of Unbiasedness and the Existence of
L-unbiased Estimators.}
\end{enumerate}

\section{Erich L. Lehmann's Bibliography}
\def\bibname{}
%



\section{Books and their translations}

\section{Epilogue}
Erich's sensitivity toward others, contagious zest for life, gentle
spirit, fundamental contributions to statistics and remarkable
contributions to human resources development, have been recorded,
chronicled and honored through various mechanisms. After his death,
Erich's life was celebrated with a memorial service that took place at
the Berkeley women's faculty club on November 9th, 2009. The
service was well attended. His family, friends, students, collaborators
and colleagues paid homage. Peter Bickel organized a memorial session
during the 2010 Joint Statistical Meetings in Vancouver (Persi
Diaconis, Juliet Shaffer and Peter Bickel speakers). The session was
very well attended with standing room only. The respect and
appreciation for Erich was international. Willem van Zwet organized a
memorial session during the 73rd IMS annual meeting in Gothenburg,
Sweden in 2010 (David Cox, Kjell Doksum, Willem van Zwet, speakers),
and Peter Bickel gave a lecture during the Latin American Congress of
Probability and Mathematical Statistics (CLAPEM) in Venezuela, November
2009, in remembrance of Erich Lehmann.

Recordings of various Erich talks are freely accessible to the public
for viewing. These include lectures he gave during the second and third
Lehmann Symposia at Rice University. Obituaries by Peter Bickel (2009) and
David Brillinger (2010) provide additional information about the life and work
of Erich L. Lehmann. Other sources that present fascinating accounts of
Erich's work and life\vadjust{\goodbreak} include Lehmann (2008B), DeGroot (1986) and Reid
(1982). A collection of selected works edited by the author will soon
be published by Springer. The \textit{Selected Works of E. L.
Lehmann} provides an extended bibliography and, through invited
vignettes, examines more closely the various facets of his work.

Stigler (2009), Rojo and Perez-Abreu (2004), Rojo (2006, 2009a, 2009b) and van
Zwet (2011) provide additional anecdotes and commentaries.

\def\bibname{References}

\printaddresses

\end{document}